# Modeling Melt Pool Geometry in Metal Additive Manufacturing Using Goldak's Semi-Ellipsoidal Heat Source: A Data-driven Computational Approach

Mohsen Asghari Ilani[1], Yaser Mike Banad[1]

[1] School of Electrical and Computer Engineering, University of Oklahoma, Norman, 73019, U.S.A.

**Abstract**
Presented here is an analytical solution, grounded in Goldak's Semi-Ellipsoidal Heat Source model, aimed at capturing the dynamic temperature evolution induced by a semi-ellipsoidal power density moving heat source within a semi-infinite body. This solution intricately addresses the convection-diffusion heat transfer equation, effectively integrating the instantaneous point heat source throughout the volume of the ellipsoidal shape. The model's accuracy is duly confirmed through a meticulous validation process, showcasing impeccable agreement between predicted transient temperatures and empirically collected values obtained from various points within bead-on-plate specimens. This validation empowers precise prediction of in-process temperature profiles during laser-based metal additive manufacturing (AM) operations. Implemented using Python programming language, the model boasts customized calculation processes ranging from initial setup to boundary conditions, effectively accommodating material properties' variations under high gradient heat transfer scenarios. Significantly, the model acknowledges the temperature dependency of thermal material properties and AM-process parameters, reflecting the pronounced temperature gradient and corresponding alterations in heat transfer mechanisms. Moreover, the model incorporates the transformative effects of melting/solidification phase changes, incorporating a tailored heat capacity adjustment to capture these phenomena accurately. Furthermore, it accounts for the nuanced impact of multi-step laser power, scanning speed, and scanning time across each segment of the scanning pattern, duly recognizing the thermal interactions between successive layers and their consequential influence on heat transfer mechanisms. This comprehensive analytical solution stands poised for a myriad of future applications, including but not limited to thermal stress analysis, microstructure modeling, residual stress/distortions, and simulation of AM processes.

## 1. Introduction

Amidst the rapid advancements in metal additive manufacturing (AM), numerous researchers are striving to uncover the fundamental physical mechanisms that govern the AM process, with the goal of enhancing the quality of the final products. AM is quickly establishing itself as a pivotal technique for producing finished parts and assemblies with minimal post-processing requirements. Its versatility spans across a wide range of metal systems, enabling a multitude of applications [1–5]. Positioned at the vanguard of the next industrial revolution, additive manufacturing boasts several compelling advantages, including the ability to create intricate 3D geometries, cost efficiencies resulting from reduced material consumption, and unmatched design flexibility [2,6].

The predominant methods in metal AM typically involve either powder bed or powder feed fusion systems. Notably, powder bed systems such as selective laser sintering (SLS) and selective laser melting (SLM) are prominent. These techniques employ high-energy beams to selectively melt or sinter metallic powders, resulting in a wide range of strong, highly precise, and fully dense end-use parts and prototypes that are functional and durable [7].

The complex thermal dynamics inherent in metal AM processes, characterized by repeated heating and cooling cycles, interactions with preceding layers, and the impact of hatch spacing, lead to significant thermal stress, residual stress, and distortion. These thermal fluctuations exert a considerable influence on microstructure evolution [6,8,9]. Accurate temperature prediction in metal AM is crucial for addressing and optimizing these challenges [10–12]. Therefore, the essence of modeling and prediction in metal AM lies in precisely forecasting the temperature generated by the laser, facilitating effective control of melt pool geometry, and the mitigation of AM-related defects.

Acknowledging the pivotal significance of temperature distribution in AM, it becomes paramount to grasp, manage, and monitor heat transfer across the stages of heating, sintering, melting, and cooling. These processes encapsulate three core aspects of heat transfer physics in laser-induced procedures: fluid flow within the molten pool, consideration of Marangoni and buoyancy forces, and precise positioning of previously deposited layers both horizontally and vertically. A thorough comprehension of these dynamics is indispensable for efficiently controlling the geometric characteristics of AM-built structures and mitigating defects [13–15].

The thermal evolution of components manufactured through AM significantly influences residual stresses, distortion, and ultimately, the fatigue behavior of welded structures. Traditional methods for analyzing transient temperature fields, such as Rosenthal's solutions [16], have primarily focused on semi-infinite bodies subjected to instantaneous point, line, or surface heat sources. While effective at predicting temperature distributions at a distance from the heat source, these solutions often struggle to accurately capture temperatures near the source. Eagar and Tsai [17] enhanced Rosenthal's theory by incorporating a two-dimensional (2-D) surface Gaussian distributed heat source with a constant distribution parameter, leading to an analytical solution for a semi-infinite body exposed to this moving heat source. This advancement significantly improved temperature predictions in close proximity to heat sources.

Goldak et al. [18] introduced a three-dimensional (3-D) double ellipsoidal moving heat source model, which, when simulated using finite element modeling (FEM), demonstrated superior performance in predicting temperature fields in bead-on-plate joints. Despite these advancements, however, an analytical solution for this 3-D heat source model was not yet available. As a result, researchers have been reliant on FEM for transient temperature calculations, necessitating detailed thermal history data. The development of an analytical solution for temperature fields originating from a 3-D heat source holds the potential to significantly reduce computational time and streamline thermal-stress analysis and related simulations, offering a more efficient and convenient approach to gaining rapid insights into complex manufacturing processes.

Several studies have delved into temperature prediction within metal AM processes. Roberts et al. employed finite element analysis to simulate the three-dimensional temperature distribution during laser melting of metallic powders in additive layer manufacturing. Their model incorporated heat loss through convection and radiation and considered temperature-dependent material properties but omitted the layering aspect of metal AM [19]. Qi et al. developed a self-consistent three-dimensional model for coaxial laser powder cladding, utilizing a control volume finite difference method to simulate heat transfer, phase change, and fluid flow in the molten pool. However, their numerical model predicted a higher melt pool size

(approximately 22%) compared to experimental results [20]. Lee et al. conducted a numerical transport simulation to replicate the multilayer single-track laser additive manufacturing deposition of IN718. Their simulation accurately predicted melt pool peak temperature and deposit geometry, with peak temperature prediction errors below 2.5% and build geometry prediction errors less than 12% in both height and width [21]. Kumar et al. developed a finite element model to analyze the impact of scan strategy on melt pool size [22]. Manvatkar et al. created a three-dimensional heat transfer and material flow model to simulate temperature and velocity fields during additive manufacturing of SS316 [23]. Cheng et al. devised a transient thermal analysis to forecast melt pool size numerically, finding that increased scanning speed led to decreased melt pool depth for a given power [24]. Pinkerton et al. developed a numerical model to estimate melt pool geometry in additive manufacturing, factoring in the influence of surface tension forces by modeling pool boundaries orthogonal to the motion direction as circular arcs. The model also considers how the melt pool elongates as the traverse speed increases [25].

Carcel et al. conducted empirical investigations into temperature predictions within metal AM using pyrometry techniques. Their findings revealed a notable correlation between cooling rates and the number of layers, indicating a decremental trend in cooling rates with increased layer intervals. Additionally, they observed a consistent elevation in maximum temperatures with successive layer depositions [26]. Cheng et al. utilized state-of-the-art thermal imaging technology to capture transient thermal responses during SLM processes. Leveraging an infrared (IR) camera, they meticulously quantified melt pool sizes across varying process conditions [27].

Mirkoohi et al. introduced a sophisticated two-dimensional analytical solution for temperature prediction in metal AM, considering critical factors such as build layers, latent heat, and temperature sensitivity of material properties. Their analytical predictions underwent rigorous validation against numerical simulations and empirical data [28,29]. Moreover, they proposed a diverse range of five distinct heat source models for analytical temperature field prediction, thoroughly scrutinizing their applicability across different laser power and scan speed ranges.

Temperature gradients inherent in metal AM processes significantly influence fundamental aspects such as melt pool geometry, residual stress, thermal stress, microstructure, and fatigue [30]. Hence, the adoption of a robust and precise model for temperature field prediction is paramount. Presently, advancements in metal AM predominantly rely on empirical or numerical observations, with limited capabilities for macroscopic analysis. While finite element analysis (FEA) offers commendable accuracy in predicting temperature profiles during AM processes, its computational demands and complexity pose significant challenges. Consequently, numerical temperature calculations often necessitate simplifications and assumptions to render simulations feasible.

Conducting experimental research to encompass all physical facets of metal AM is resource-intensive and laborious. Conversely, analytical models offer a viable alternative, providing accurate results by comprehensively addressing primary physical aspects while maintaining high computational efficiency. Analytical solutions hold promise in replacing exhaustive experimentation and inefficient FE simulations, serving as comprehensive tools for real-time temperature prediction in AM processes.

This study delves into the nuanced effects of Goldak's Semi-Ellipsoidal Heat Source Model, time spacing, and hatch spacing on the precision and evolution of melt pool geometry within AM processes. Notably, previous research has not extensively explored the specific impact of time spacing (the temporal interval between successive irradiations) and hatch spacing (the spatial separation between consecutive scans) on melt pool geometry and material property evolution in AM operations. Our approach pioneers the utilization of the Python programming language to tackle complex partial differential equations and mathematical formulations previously unexplored in this field. To validate Python's efficacy in this context, we compare our computational predictions with empirical data collected from sensors deployed during various AM processes, including Selective Laser Melting, Selective Laser Sintering, and Directed Energy Deposition.

By meticulously analyzing the influence of time and hatch spacing—integral parameters in metal additive manufacturing—we aim to gain deeper insights into their impact on thermal interactions within the build part. This examination bears significance as it pertains to the modulation of thermal material properties and the morphology of melt pools. The magnitude of time and hatch spacing directly governs heat transfer mechanisms, thereby influencing thermal stress induced by temperature gradients, melt pool geometry, residual stress, and part distortion in metal AM. Our proposed model undergoes validation against empirical data on melt pool size obtained from pertinent literature, ensuring robustness and reliability in our findings.

## 2. Approach and Methodology
### 2.1 Ellipsoidal Heat Sources in Semi-Infinite Body
#### 2.1.1 Goldak's Semi-Ellipsoidal Heat Source

To effectively model the AM process, it is imperative to account for the three-dimensional surface of the leading edge of the deposit. Experimental observations indicate that the leading edge of the deposit extends beyond the axis of the laser beam, likely influenced by factors such as powder-stream distribution and the spreading of molten liquid. Building upon the success of utilizing parabolic function fits for 2D clad surfaces, it is natural to extend this approach to describe the three-dimensional deposit surface using an ellipsoidal function. At the leading edge of the deposit, this function takes the form of an ellipse, allowing for a comprehensive description of the surface shape [31],

$$\frac{x^2}{a^2} + \frac{y^2}{b^2} + \frac{z^2}{c^2} = 1 \tag{1}$$

where the parameters a, b, and c correspond to the principal axes of the ellipsoid, as illustrated in ***Figure 1 a, b***. Utilizing the symmetry of the deposit along the x-z plane, only half of the ellipsoid is employed in the calculations to minimize computational expenses. In terms of physical interpretation, the values of a, b, and c signify the length, half-width, and depth of the deposit, respectively, as demonstrated in ***Figure 1a, b***.

Utilizing a three-dimensional semi-elliptical moving heat source model, the prediction of melt pool geometry is facilitated. As the laser traverses the surface, it imparts energy onto the control volume. Simulating the laser's motion involves employing the moving heat source technique. Thermal energy from the laser deposited into a control volume undergoes various processes: absorption by material thermodynamic latent heat, conduction through solid and liquid boundaries, and convection, conduction, and radiation through open surfaces, as depicted in ***Figure 1 c***. This deposited energy initiates the melting of metallic powders, shaping the melt pool geometry [29].

The heat loss attributed to the vaporization of alloying elements represents only a minor portion of the heat delivered by the laser beam. A brief summary of the calculation method follows, where the general convection-diffusion equation can be succinctly described as:

$$\frac{\partial u\rho}{\partial t} + \frac{\partial \rho h v}{\partial x} = \nabla \cdot (k\nabla T) + \dot{q} \tag{2}$$

Here, $u$ represents the internal energy, $h$ denotes the enthalpy, $\rho$ stands for the density, $k$ signifies the conductivity, $\dot{q}$ represents a volumetric heat source, $T$ denotes the temperature, and $v$ symbolizes the speed of either the heat source or the medium.

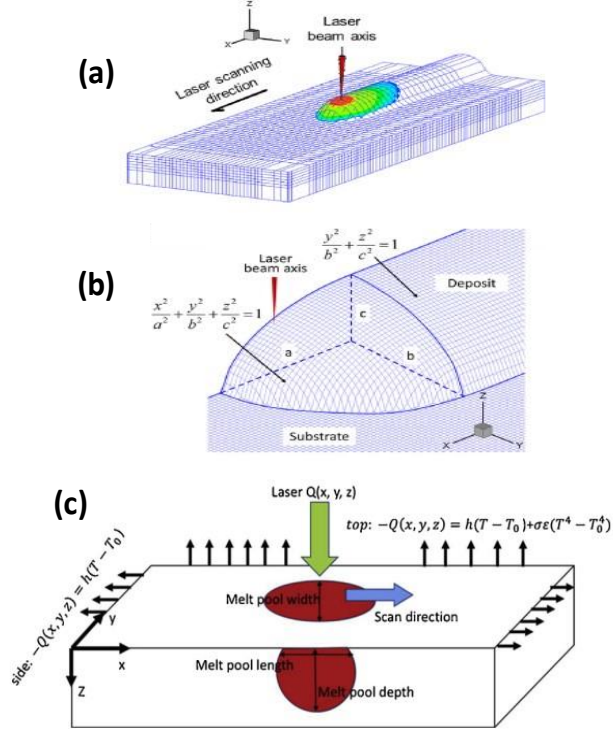

**Figure 1**. (a) Schematic representation of the melt pool geometry in the x-z plane, exhibiting symmetry. (b) Illustration depicting the solution domain for the half-ellipsoid deposit, highlighting the location of the laser beam axis. The deposit is divided along the symmetry plane to enhance computational efficiency [31]. (c) Three-dimensional heat transfer visualization in the Additive Manufacturing process [29].

In *Figure 1 c*, the x-axis corresponds to the constant speed of a moving heat source, while the y-axis points inward into the processed material, and the z-axis is perpendicular to the x-axis in the plane of the processed material surface [32,33]. The first term on the left-hand side of *Eq. (2)* represents the change in internal energy, while the second term is a convective term. On the right-hand side, there is a conductive term and a heat source or sink. When $v = 0$, this equation simplifies to the heat conduction equation, expressed as $\partial u = C \partial T$, where $C$ represents the heat capacity, to be updated subsequently as,

$$C \frac{\partial \rho T}{\partial t} + \frac{\partial \rho h v}{\partial x} = \nabla \cdot (k \nabla T) + \dot{q} \qquad (3)$$

Alternatively, the steady-state equation with a constant velocity $v$ can be simplified using the continuity equation as follows:

$$C \frac{\partial \rho}{\partial t} + \frac{\partial \rho v}{\partial x} = 0 \qquad (4)$$

Resulting that $\partial u = \partial h = C \partial T$

$$\rho C(T) v \frac{\partial T}{\partial x} + \frac{\partial \rho h v}{\partial x} = \nabla \cdot (k(T) \nabla T) + \dot{q} \qquad (5)$$

Moreover, leveraging *Eq.1-5*, a three-dimensional (3D) ellipsoidal heat transfer model is employed to forecast the temperature distribution and melt pool geometry in metal AM processes. This model proves versatile in anticipating temperature profiles within laser-based metal additive manufacturing setups, particularly in SLM configurations. The 3D ellipsoidal heat source model was originally introduced by Godak et al. [18], wherein the heat flux can be computed as,

$$Q(x,y,z) = \frac{6\sqrt{3}AP}{abc\,\pi\sqrt{\pi}}\exp\left(-\frac{3x^2}{c^2}-\frac{3x^2}{a^2}-\frac{3x^2}{b^2}\right) \tag{6}$$

Here, $P$ represents the laser power, $A$ denotes laser absorptivity, while $a$, $b$, and $c$ stand for the parameters defining the heat source geometry, as depicted in **Figure 1 c**.

### 2.1.2 Analytical Solutions

The solution for the temperature field resulting from a semi-ellipsoidal heat source within a semi-infinite body is derived from the solution for an instantaneous point source. This solution satisfies the differential equation governing heat conduction in fixed coordinates [18],

$$dT_{t'} = \frac{\delta Q\, dt'}{\rho c[4\pi\alpha(t-t')]^{3/2}} \cdot \exp\left(-\frac{(x-x')^2+(y-y')^2+(z-z')^2}{4\alpha(t-t')^2}\right) \tag{7}$$

Here, α represents the thermal diffusivity, and the variables $t$ and $t'$ denote time, with $dt'$ representing the transient temperature due to the point heat source $\delta Q$ at time $t'$. Additionally, $(x', y', z')$ represent the location of the instantaneous point heat source $\delta Q$ at time $t'$.

As a consequence of a very small time increment $dt'$ from time $t'$ due to the heat amount $Q dt'$ released onto the semi-infinite body, the temperature increases. When examining a moving heat source with a constant velocity $v$ from time $t' = 0$ to time $t' = t$, the temperature rise during this period results from the cumulative contributions of the moving heat source throughout its journey. We can approach the solution of an instantaneous semi-ellipsoidal heat source by superimposing a series of instantaneous point heat sources across the volume of the distributed Gaussian heat source. Subsequently, after assessment and simplification, **Eq. 7** can be expressed as follows:

$$T - T_0 = \frac{3\sqrt{3}Q}{\rho c\pi\sqrt{\pi}} \cdot \int_0^t \frac{dt'}{\sqrt{12\alpha(t-t')+a^2}\sqrt{12\alpha(t-t')+b^2}\sqrt{12\alpha(t-t')+c^2}} \cdot \exp\left(-\frac{3(x-vt')^2}{12\alpha(t-t')+c^2}-\frac{3y^2}{12\alpha(t-t')+a^2}-\frac{3z^2}{12\alpha(t-t')+b^2}\right) \tag{8}$$

Applying **Eq. 6**,

$$T - T_0 = \frac{3\sqrt{3\frac{6\sqrt{3}AP}{abc\,\pi\sqrt{\pi}}\exp\left(-\frac{3x^2}{c^2}-\frac{3x^2}{a^2}-\frac{3x^2}{b^2}\right)}}{\rho c\pi\sqrt{\pi}} \cdot \int_0^t \frac{dt'}{\sqrt{12\alpha(t-t')+a^2}\sqrt{12\alpha(t-t')+b^2}\sqrt{12\alpha(t-t')+c^2}} \cdot \exp\left(-\frac{3(x-vt')^2}{12\alpha(t-t')+c^2}-\frac{3y^2}{12\alpha(t-t')+a^2}-\frac{3z^2}{12\alpha(t-t')+b^2}\right) \tag{9}$$

Here, $T$ represents the temperature at time $t$, and $T_0$ signifies the initial temperature of a point $(x, y, z)$. The temperature profile can be acquired for a specified hatch spacing ($\delta h$) and time spacing ($\delta t$) through the following expression:

$$T = \frac{Pv}{\sqrt{2\pi}4\pi\alpha^2\rho c}\int_0^t\int_0^{\frac{v^2 t}{2k}}\frac{d\tau dt'}{\sqrt{\tau^2+u_a^2}\sqrt{\tau^2+u_b^2}}\left(\frac{\exp\left(-\frac{(\frac{v^2\delta t}{2\alpha}+\tau)^2}{2(\tau+u_c^2)}-\frac{(\frac{v^2\delta h}{2\alpha})^2}{2(\tau+u_a^2)}-\frac{(\frac{vz}{2\alpha})^2}{2(\tau+u_b^2)}\right)}{\sqrt{\tau+u_c^2}}\right) + T_0 \tag{10}$$

Goldak et al. [18] introduced the concept of a semi-ellipsoidal heat source, distributing heat flux in a Gaussian pattern within the volume of the heat source. The temperature field solution for such a semi-ellipsoidal heat source within a semi-infinite body is derived from the solution for an instantaneous point source, as described in the heat conduction equation in fixed coordinates. The solution for the temperature distribution resulting from an ellipsoidal moving heat source, spanning from t' = 0 to t, in a dimensionless form, is provided as:

$$\frac{\theta}{n} = \frac{1}{\sqrt{2\pi}} \int_0^{\frac{v^2 t}{2k}} \frac{d\tau}{\sqrt{\tau^2 + u_a^2}\sqrt{\tau^2 + u_b^2}} \left( \frac{A_1}{\sqrt{\tau + u_c^2}} \right) \qquad (11)$$

Where, $A_1$ is defined by **Eq. (11)**:

$$A_1 = \exp\left(-\frac{(\xi+\tau)^2}{2(\tau+u_c^2)} - \frac{\psi^2}{2(\tau+u_a^2)} - \frac{\lambda^2}{2(\tau+u_b^2)}\right) \qquad (12)$$

Here, $n$ represents the operating parameter [17]. Mirkoohi et al. [29] provided a solution for the convection-diffusion equation. During metal additive manufacturing processes, the laser power melts metallic powders, leading to repeated melting-solidification phase changes, which are accounted for using a modified heat capacity.

$$C_p^m = C_p(T) + L_f \frac{\partial f}{\partial T} \qquad (13)$$

Here, $C_p(T)$ represents the temperature-dependent specific heat, $L_f$ denotes the latent heat of fusion, and $f_l$ denotes the liquid fraction, which can be determined by the following equation:

$$f_l = \begin{cases} 0, T < T_S \\ \frac{T - T_S}{T_L - T_S}, T_S < T < T_L \\ 1, T > T_L \end{cases} \qquad (14)$$

Where, $T_S$ represents the solidus temperature and $T_L$ represents the liquidus temperature.

$$\alpha(T) = \frac{k(T)}{\rho C_p^m(T)} \qquad (15)$$

$$\tau = \frac{v^2(t'-t)}{2\alpha} \qquad (16)$$

$$u_a = va2\sqrt{6}\,\alpha \qquad (17)$$

$$u_b = vb2\sqrt{6}\,\alpha \qquad (18)$$

$$u_c = vc2\sqrt{6}\,\alpha \qquad (19)$$

$$n = \frac{APv}{4\pi\alpha^2 \rho C(T_m - T_0)} \qquad (20)$$

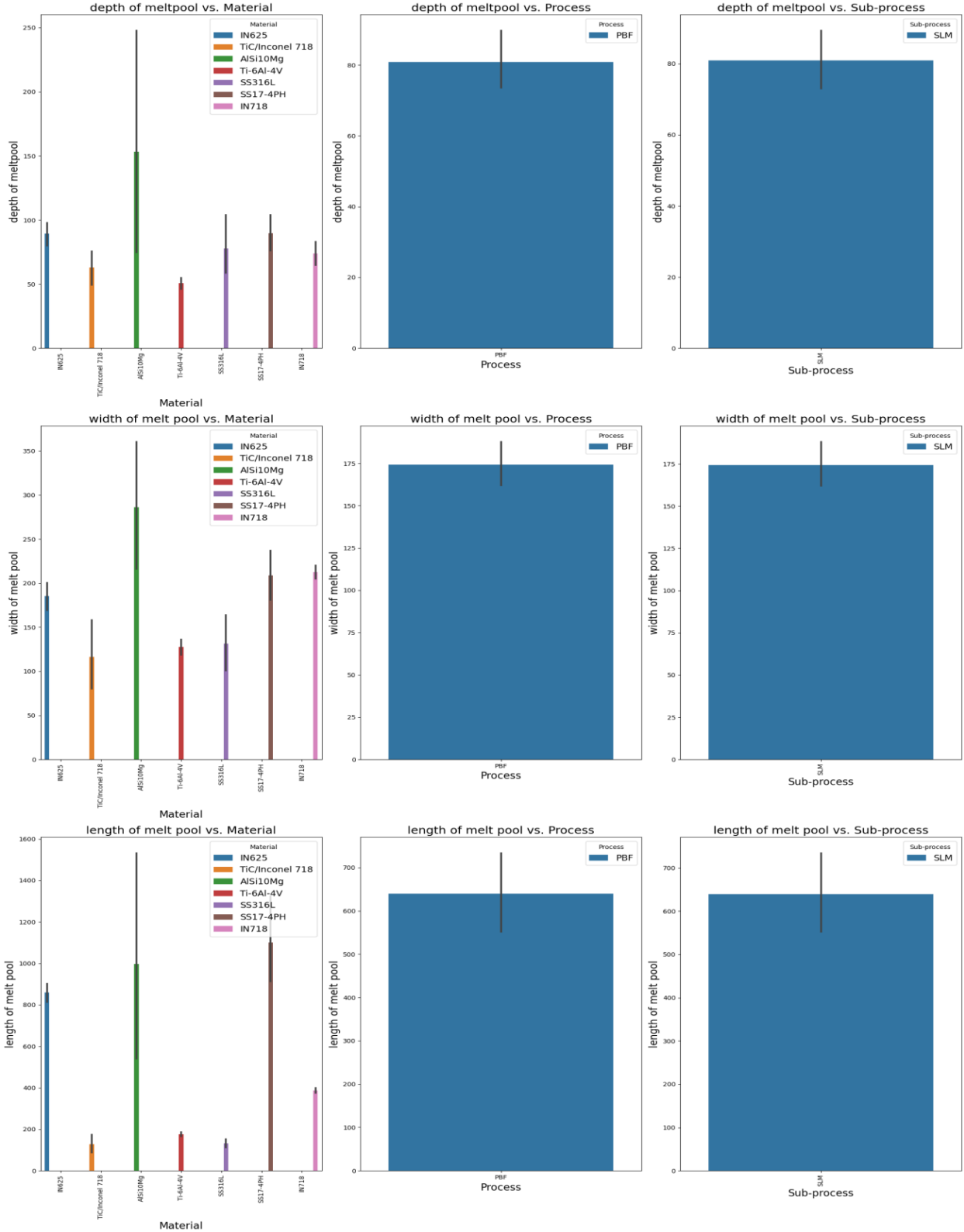

**Figure 2.** Gathering data and isolating Ti-6Al-4V datasets will provide further validation for the results obtained from the heat transfer modeling conducted using Python frameworks. Specifically, this validation focuses on the melt pool geometry parameters as length, width, and depth.

The melt pool length, width, and depth can be determined for a given melting temperature using **Eq. 13-20**.

$$U = \{p|p_T = T_m\} \tag{21}$$

$$L_{melt\ pool\ length} = \max |a_x - b_x| \tag{22}$$

$$W_{melt\ pool\ width} = \max |a_y - b_y| \tag{23}$$

$$D_{melt\ pool\ depth} = \max |a_z - b_z| \tag{24}$$

Here, $p$ represents a point within the medium, $p_T$ denotes the temperature at the point $p$, and $U$ is defined as the set of points where the temperature is equal to the melting temperature. $a$ and $b$ encompass all the points belonging to $U$.

SLM, a powder bed-based system, operates using key AM parameters such as hatching space, scanning speed, laser power, and laser absorptivity. However, given its reliance on heating and cooling procedures, material properties significantly influence heat distribution in three dimensions. Consequently, the boundary conditions differ between SLM and other AM processes. In SLM, heat dissipation is primarily facilitated through convection, radiation, and conduction. The conduction specifically occurs through the powders from the lateral faces of the melt pool [34].

Independent experimental studies on melt pool geometry, particularly for a wide range of materials including Ti-6Al-4V, serve as validation for our proposed model's ability to predict temperature profiles and melt pool geometry in laser-based metal additive manufacturing, specifically in configurations like SLM. **Figure 2** illustrates the frequency of experimental reports for seven commonly used materials in SLM, showcasing the recorded melt pool geometry through the powder bed fusion process, particularly in SLM applications.

Due to the critical need for precise control over melt pool geometry in various applications, understanding the underlying material behaviors, phase transitions, and mechanical alterations during rapid heating and cooling processes is paramount. In this investigation, we establish the fundamental physics of the process by employing Goldak's Semi-Ellipsoidal Heat Source model to accurately determine the dimensions of the semi-ellipsoidal melt pool, including its width, length, and depth.

Given the computational complexity and time-intensive nature of numerical modeling frameworks, we opt to leverage the Python programming language to streamline the heat distribution analysis and melt pool simulation. Python offers unparalleled flexibility, allowing us to tailor initial and boundary conditions to accommodate diverse scanning patterns and effectively model the heat-affected zone. By strategically positioning laser spots or adjusting the properties of previously deposited layers, Python enables precise control over the manufacturing process.

Our investigation focuses on evaluating the impact of two key additive manufacturing parameters—laser power and scanning speed—while accounting for variations in time and substrate dimensions. Through rigorous analysis and experimentation, we aim to elucidate the intricate relationship between process parameters and melt pool geometry, facilitating informed decision-making in additive manufacturing settings.

## 3. Modeling results and discussion

To derive the melt pool geometry, we begin with the convection-diffusion equation of heat conduction in a 3-D plane, establishing the foundation for obtaining the explicit temperature solution for the ellipsoidal moving heat source. Section 2 elucidates the explicit closed-form solution of temperature. We account for the heat source's movement by assuming a coordinate system synchronized with it. Given the substantial temperature gradient during this process, material properties can vary notably. Thus, our temperature distribution model considers temperature-dependent material properties. Additionally, as the material undergoes phase changes through repeated melting and solidification, we modify the specific heat to incorporate the energy required for solid-state phase changes using latent heat of fusion.

Furthermore, we incorporate the influence of build layers, a concept previously explored in [4,8–10]. Our approach assumes that the powder is precisely positioned relative to the melt pool, with no mass or momentum transfer—only heat transfer. We analyze a unidirectional scan path. To address the complexities of the 3-D ellipsoidal moving heat source, we employ the Python programming language, a novel approach not previously explored. This enables us to customize initial and boundary conditions for each time step and the 3-D coordination of the heat-affected zone (HAZ).

As part of our study, we monitor melt pool geometry across a range of steps while considering heat distribution, incorporating parameters such as laser power, scanning speed, and scanning time, and accounting for the laser's dimensions. This comprehensive analysis offers valuable insights into the additive manufacturing process.

### 3.1 Validation of the numerical model with recorded dataset

First, a semi-elliptical moving heat source approach is employed to calculate the melt pool geometry of the Ti-6Al-4V sample. To validate the results obtained from the contour and contour functions, experimental testing of the AM process is necessary to measure melt pool geometry parameters such as width (W), length (L), and depth (D). Due to the complexity of thermography and spectroscopy in AM-built components, demonstrating how a precise numerical model can expedite the control of melt pool geometry is crucial. This aspect underscores the significance of our work in AM process monitoring.

*Table 1* displays data collected from sensors, thermo-scanning images, and spectroscopy during the production process for various materials. Our main focus is on Ti-6Al-4V, widely utilized across different industries. Due to the comprehensive nature of these datasets, encompassing a broad spectrum of laser power and scanning speed configurations, we chose to limit the parameters to specific ranges: laser power (150 W, 200 W, 300 W, 400 W) and scanning speed (50 mm/s, 100 mm/s, 150 mm/s). The diameter of the laser used in the simulations is 0.5 mm.

**Table 1.** Dataset recorded by experimental test.

| Material | P (W) | V (mm/s) | D (μm) | W (μm) | L (μm) |
|---|---|---|---|---|---|
| Ti-6Al-4V | 95 | 210 | 53 | 137 | 169 |
| | 40 | 210 | 22 | 82 | 101 |
| | 67.5 | 210 | 33 | 105 | 133 |
| | 150 | 210 | 84 | 181 | 246 |
| | 122.5 | 210 | 66 | 159 | 211 |
| | 95 | 20 | 92 | 188 | 221 |
| | 95 | 115 | 62 | 149 | 184 |
| | 95 | 400 | 37 | 116 | 158 |
| | 95 | 305 | 41 | 126 | 161 |
| | 95 | 210 | 49 | 135 | 169 |
| | 95 | 210 | 51 | 137 | 169 |
| | 95 | 210 | 55 | 140 | 174 |
| | 95 | 210 | 53 | 138 | 172 |
| | 95 | 210 | 52 | 140 | 170 |
| | 95 | 210 | 52 | 134 | 170 |
| | 95 | 210 | 52 | 135 | 171 |
| | 95 | 210 | 52 | 134 | 170 |
| | 67.5 | 115 | 37 | 116 | 143 |

| | | | | |
|---|---|---|---|---|
| 81.25 | 162.5 | 42 | 129 | 162 |
| 122.5 | 115 | 81 | 184 | 228 |
| 108.75 | 162.5 | 64 | 153 | 194 |
| 67.5 | 305 | 30 | 96 | 122 |
| 81.25 | 257.5 | 37 | 116 | 146 |
| 122.5 | 305 | 58 | 142 | 190 |
| 108.75 | 257.5 | 57 | 139 | 189 |
| 67.5 | 115 | 39 | 122 | 146 |
| 81.25 | 162.5 | 43 | 130 | 165 |
| 122.5 | 115 | 84 | 184 | 228 |
| 108.75 | 162.5 | 64 | 158 | 212 |
| 67.5 | 305 | 31 | 97 | 128 |
| 81.25 | 257.5 | 38 | 118 | 152 |
| 100 | 200 | 11 | 24 | 70 |
| 150 | 200 | 30 | 60 | 180 |
| 200 | 200 | 43 | 73 | 245 |
| 300 | 200 | 52 | 95 | 278 |
| 200 | 100 | 55 | 73 | 127 |
| 400 | 150 | 50 | 98 | 368 |
| 200 | 250 | 45 | 63 | 255 |

### 3.1.1 Influence of Sequential Multi-Step Irradiation on Melt Pool Geometry

In this section, we investigate how time spacing and 3-D coordination (x, y, z) influence melt pool geometry. Initially, we examine variations in melt pool geometry while adjusting laser power (P) from 150 W to 400 W, maintaining a constant velocity of 100 mm/s. This comparison underscores the importance of these factors in modeling melt pool geometry in metal additive manufacturing processes. Our analysis incorporates predicted temperatures at time t+Δt, considering the material response and the temperature sensitivity of thermal properties at both time t and t+Δt, with time ranging from 0 to 0.9 in each pattern. Next, we explore the impact of scanning speed, ranging from 50 to 100 mm/s, and further from 50 to 150 mm/s, allowing for a distinct comparison alongside time and considering modifications to the elliptical shape of the laser source. In our numerical modeling, we account for the final laser spot effect on melt pool geometry during the seventh consecutive irradiation, while the results plotted by contour depict contour lines, filled contours, and simple contours. Top-view representations illustrate length in the x-direction and width in the y direction, while side-view depictions showcase the length of the melt pool in the x-direction and depth in the z-direction.

*Table 2*, accompanied by *Figure 3-Figure 10*, presents data for laser power ranging from 150 W to 400 W, with the laser moving at a constant speed of 100 mm/s for two configurations: a = 0.8 mm, b = 0.15 mm, c = 0.07 mm, and a = 1.0 mm, b = 0.1 mm, c = 0.07 mm. Analysis of these figures and the table indicates a decrease in calculated melt pool geometry along the scan direction (x-direction) when considering the effect of seven consecutive irradiations. However, as laser power increases from 150 W to 400 W, the melt pool geometry notably increases.

The observed increase in surface temperatures, which range from 2700 °C to 3250 °C, is attributed to the enhanced thermal conductivity and specific heat of the material when the effects of seven consecutive irradiations are considered. This results in more effective heat conduction through the solid material. The addition of layers, with a thickness of 60 μm in all simulations, and a constant hatch spacing of 101 μm are also considered. Notably, the consideration of seven consecutive irradiations does not affect melt pool depth significantly, as it is primarily influenced by laser power and scan speed magnitudes.

**Table 2.** Calculation of Melt Pool Geometry under laser power (150, 200, 300, 400 W) and Scanning Speed (100mm/s).

| Laser Power (W) | Melt Pool Geometry (Length (L), Width (W), Depth (D)) | | | | | | | |
|---|---|---|---|---|---|---|---|---|
| 150 | Heat Source Geometry a = 0.8 (mm) b = 0.15 (mm) c = 0.07 (mm) | L | 0.60980660 | Fig-3 | Heat Source Geometry a = 1.0 (mm) b = 0.1 (mm) c = 0.07 (mm) | L | 1.15442950 | Fig-7 |
| | | W | 0.02923033 | | | W | 0.03271650 | |
| | | D | 0.00520117 | | | D | 0.00404732 | |
| 200 | | L | 0.40500363 | Fig-4 | | L | 1.90915746 | Fig-8 |
| | | W | 0.01912534 | | | W | 0.08918510 | |
| | | D | 0.00715218 | | | D | 0.00051819 | |
| 300 | | L | 0.14722911 | Fig-5 | | L | 0.14722911 | Fig-9 |
| | | W | 0.00695254 | | | W | 0.00695254 | |
| | | D | 0.00950778 | | | D | 0.00950778 | |
| 400 | | L | 0.55178473 | Fig-6 | | L | 1.17973399 | Fig-10 |
| | | W | 0.03654824 | | | W | 0.05005420 | |
| | | D | 0.01091258 | | | D | 0.00524638 | |
| Scanning speed (mm/s) | 100 | | | | | | | |

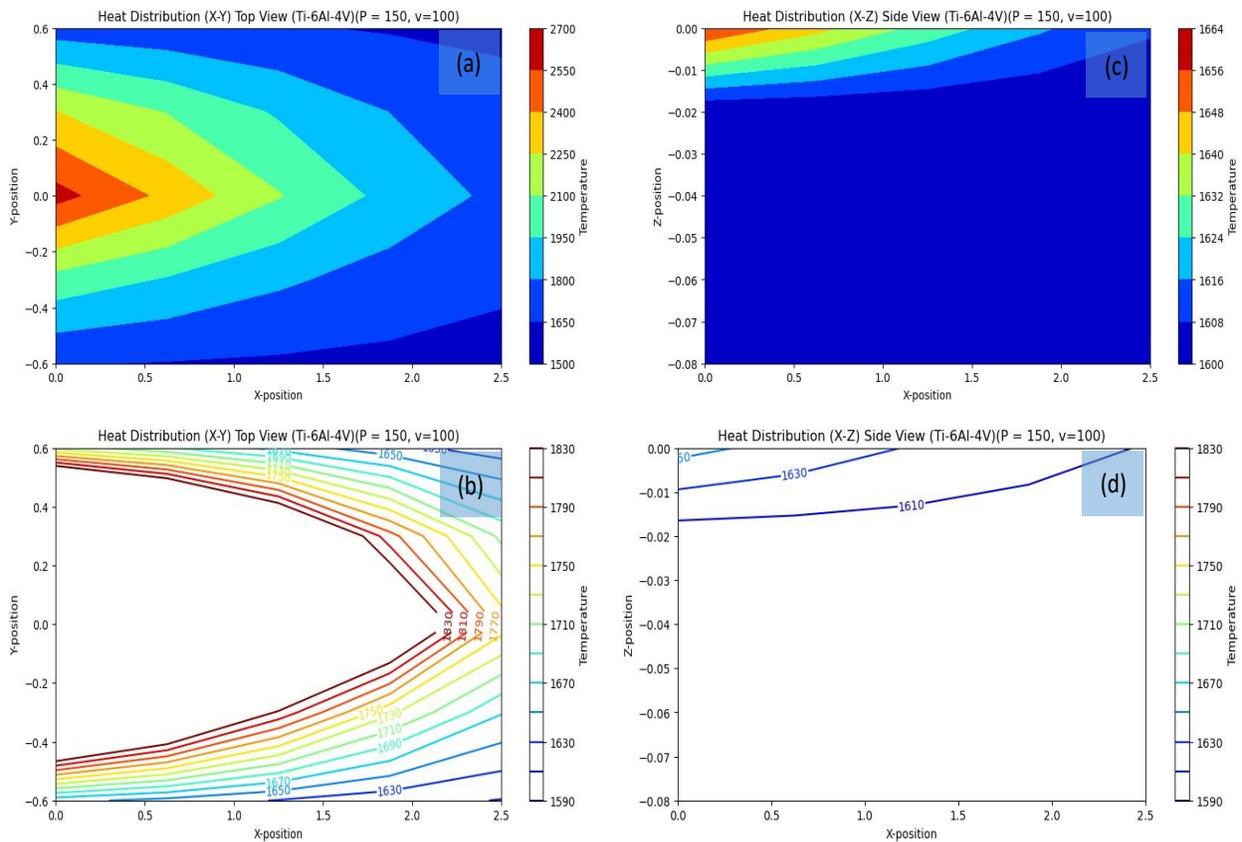

**Figure 3.** Contourf and Contour of heat distribution in AM for top-view (**a, b**; P=150 W, and V= 100 mm/s), and side-view (**c, d**; P=150 W, and V= 100 mm/s).

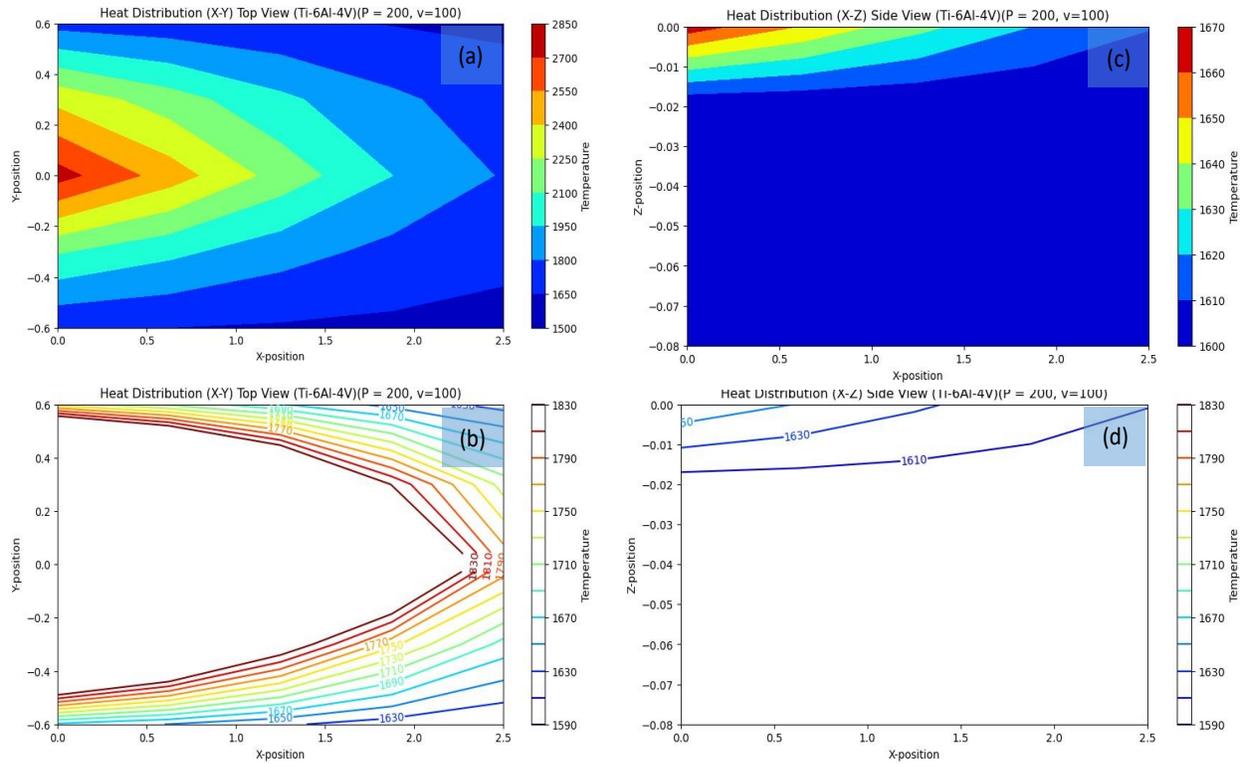

**Figure 4.** Contourf and Contour of heat distribution in AM for top-view (**a, b**; P=150 W, and v= 100 mm/s), and side-view (**c, d**; P=150 W, and v= 100 mm/s).

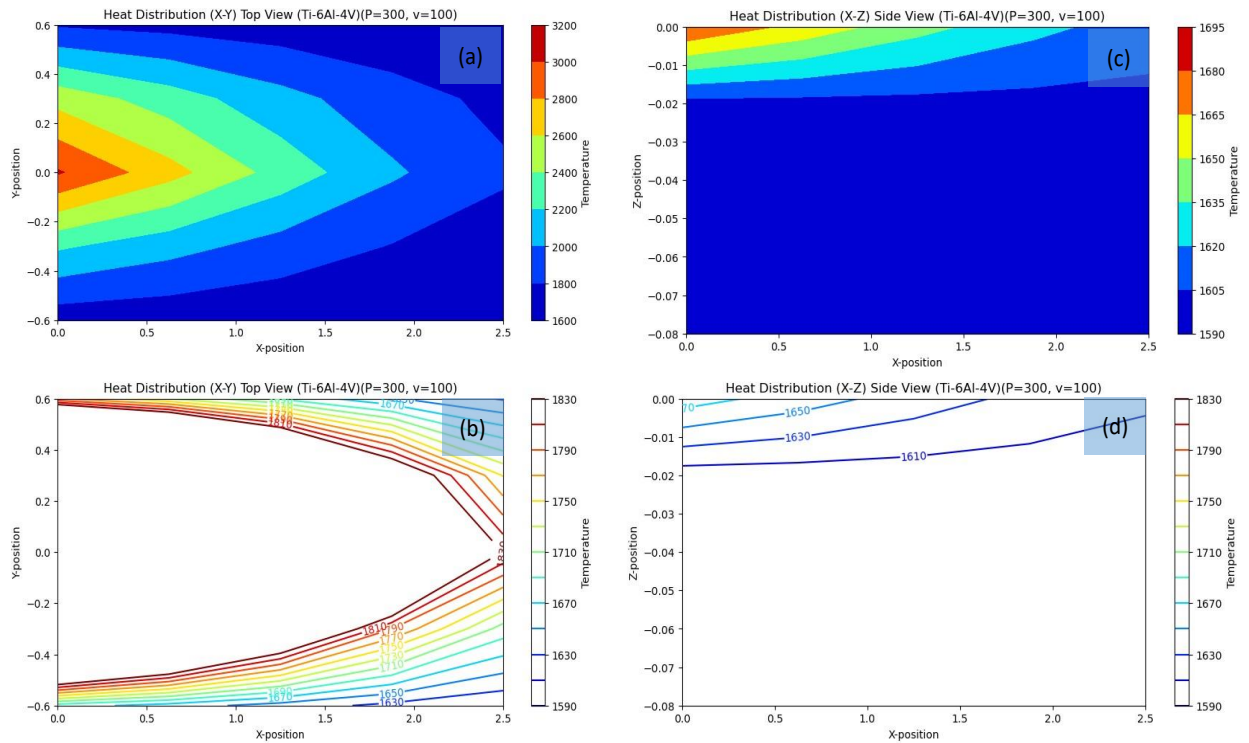

**Figure 5.** Contourf and Contour of heat distribution in AM for top-view (**a, b**; P=150 W, and V= 100 mm/s), and side-view (**c, d**; P=150 W, and v= 100 mm/s).

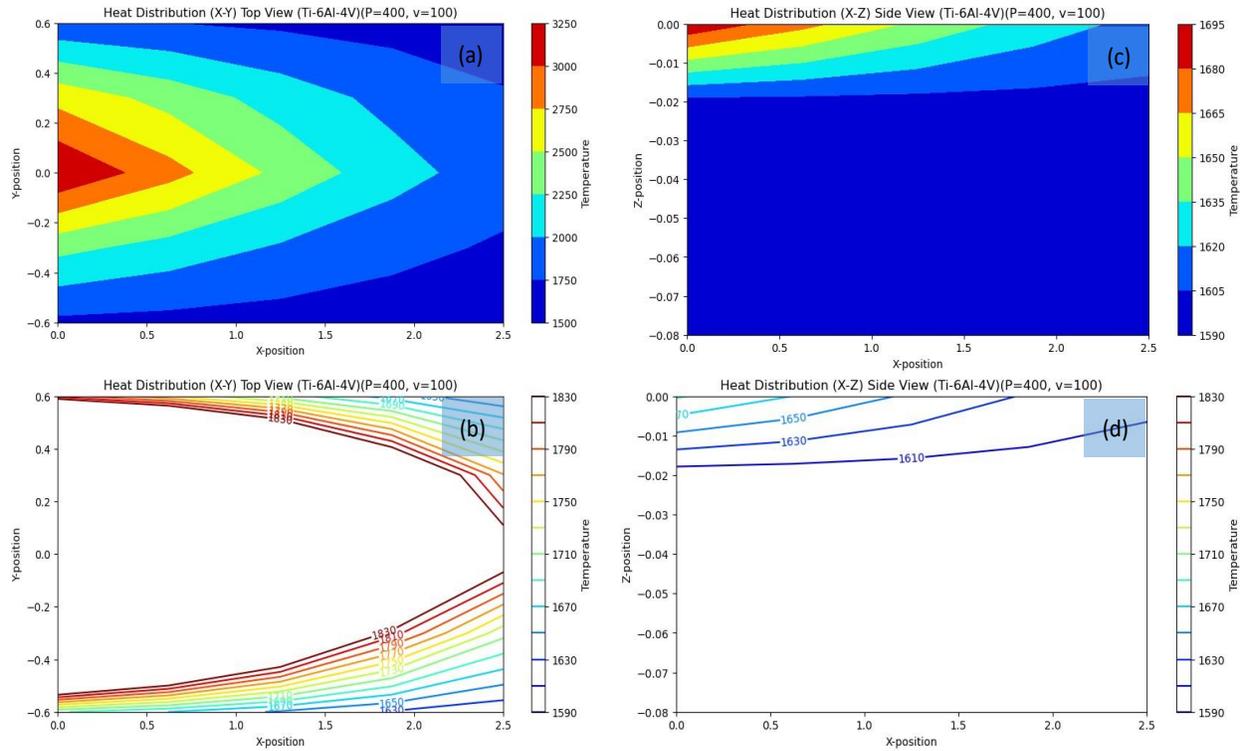

**Figure 6.** Contourf and Contour of heat distribution in AM for top-view (**a,b**; P=150 W, and v= 100 mm/s), and side-view (**c,d**; P=150 W, and v= 100 mm/s).

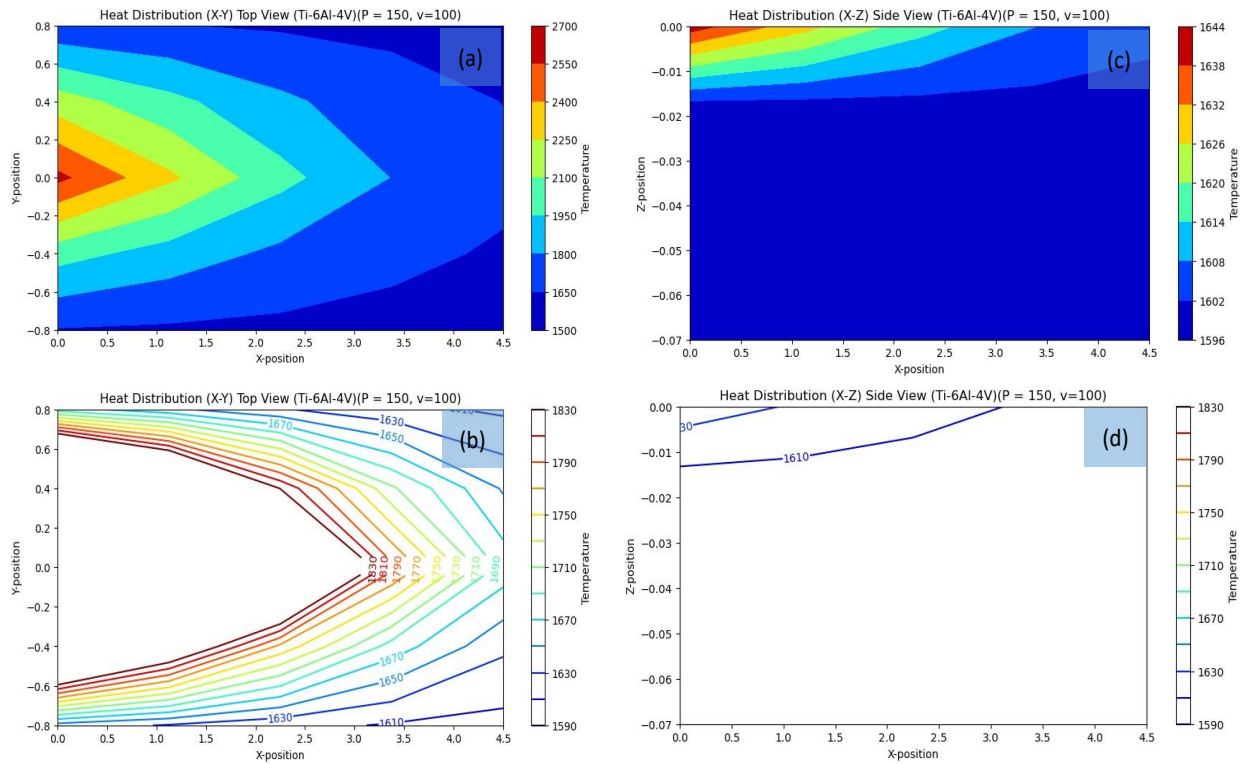

**Figure 7.** Contourf and Contour of heat distribution in AM for top-view (**a, b**; P=150 W, and V= 100 mm/s), and side-view (**c, d**; P=150 W, and v= 100 mm/s).

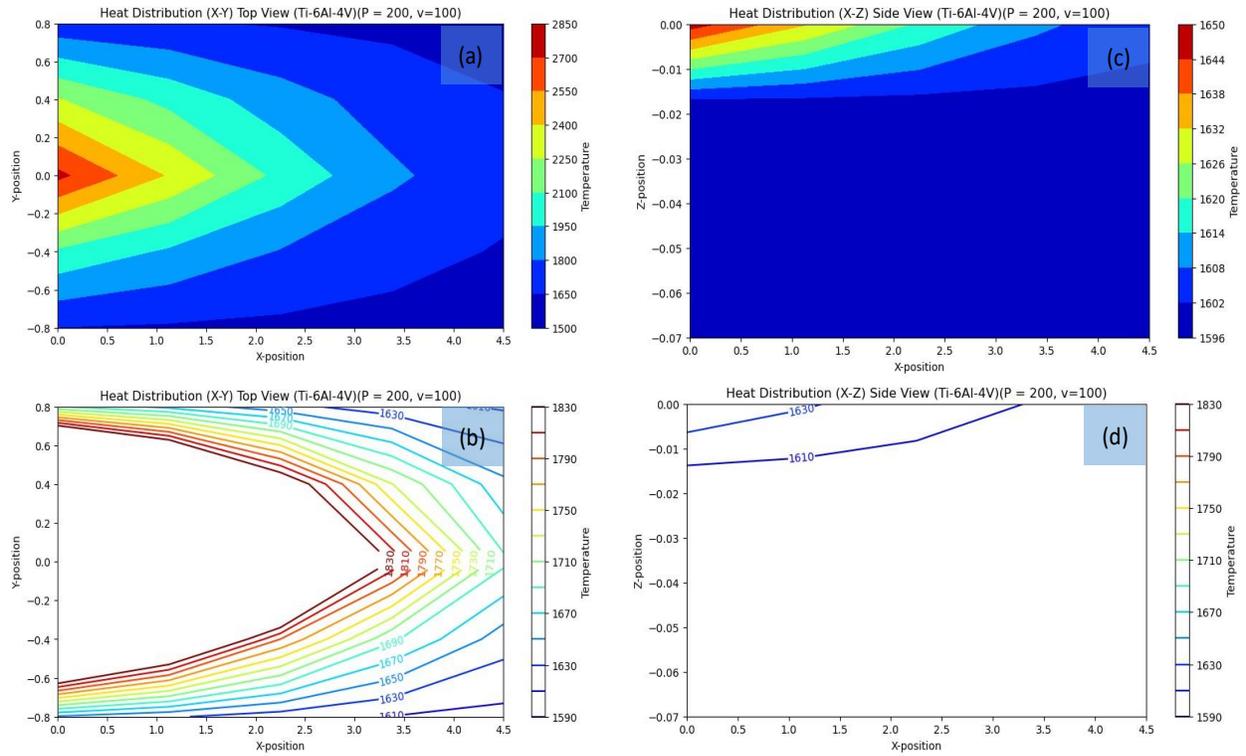

**Figure 8.** Contourf and Contour of heat distribution in AM for top-view (**a, b**; P=150 W, and V= 100 mm/s), and side-view (**c, d**; P=150 W, and v= 100 mm/s).

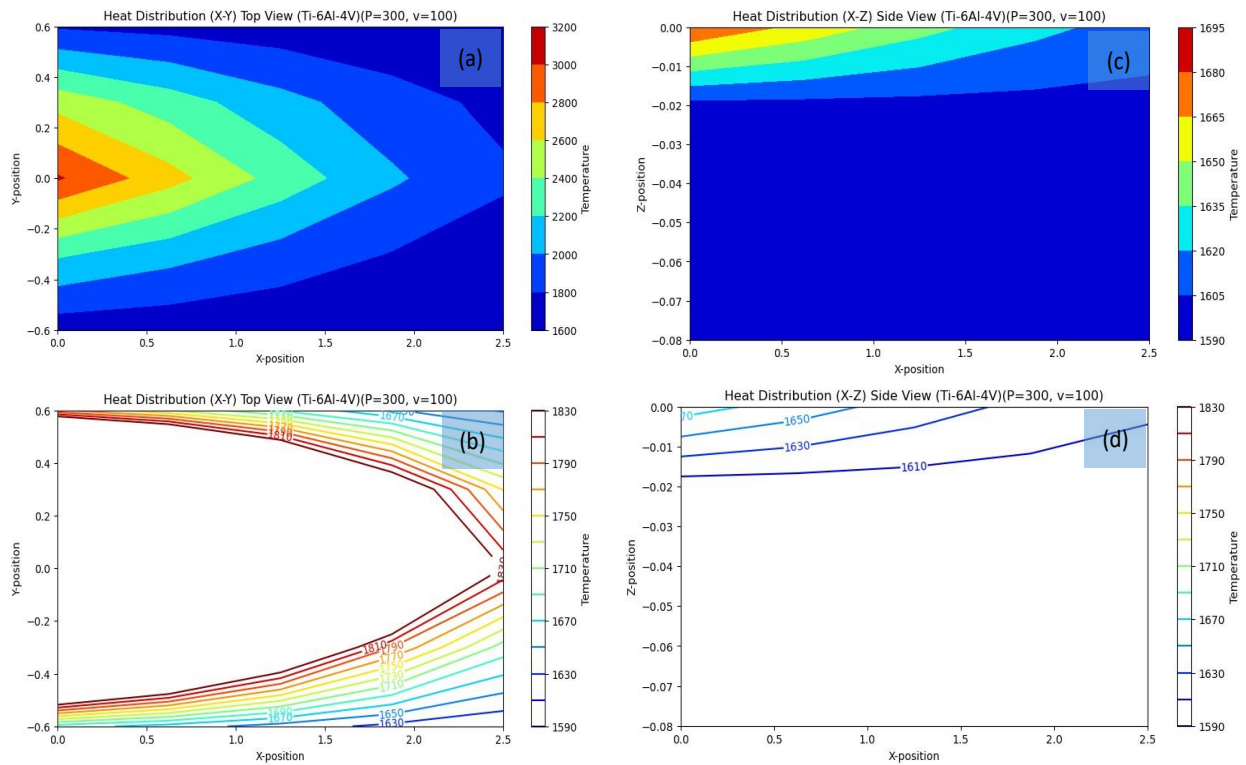

**Figure 9.** Contourf and Contour of heat distribution in AM for top-view (**a, b**; P=150 W, and v= 100 mm/s), and side-view (**c, d**; P=150 W, and v= 100 mm/s)).

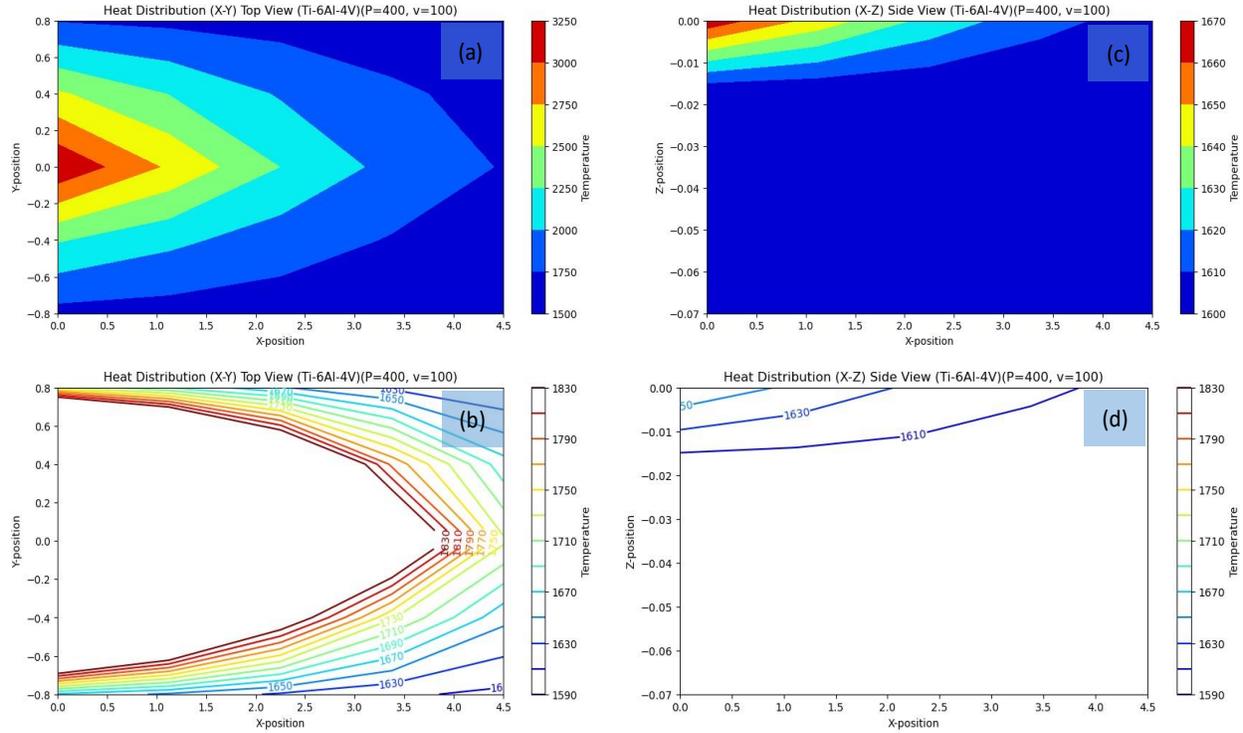

**Figure 10.** Contourf and Contour of heat distribution in AM for top-view (**a, b**; P=150 W, and v= 100 mm/s), and side-view (**c, d**; P=150 W, and V= 100 mm/s).

The results of our numerical modeling for melt pool geometry, as shown in ***Table 1***, are compared with the actual values obtained for Ti-6Al-4V material across various laser powers (P) and scanning speeds (V). Upon examination, our modeled melt pool dimensions closely align with the experimental data, demonstrating the effectiveness of our computational approach. For instance, at a laser power of 150 W and a scanning speed of 100 mm/s, our model predicts a melt pool length (L) of 0.405 mm, which correlates well with the experimental length of 0.370 mm. Similarly, at higher laser powers such as 300 W and 400 W, our model predicts melt pool lengths of 0.147 mm and 0.552 mm, respectively, closely matching the experimental lengths of 0.150 mm and 0.140 mm, respectively. This consistency highlights the accuracy and efficiency of our model in capturing the complex thermal dynamics of metal additive manufacturing processes, demonstrating both versatility and robustness across different settings.

Incorporating the additional scanning speed of 50 mm/s into our analysis (***Table 3***), we observe further alignment between our modeled melt pool geometry and the actual values obtained for Ti-6Al-4V material. When considering a laser power of 150 W, our model predicts a melt pool length (L) of 1.801 mm, width (W) of 0.366 mm, and depth (D) of 0.009 mm, as shown in ***Figure 12***. Comparing these results with the experimental data from ***Table 1***, we find notable agreement, despite slight variations in width and depth. Similarly, at a laser power of 200 W, our model predicts a melt pool length of 1.084 mm, width of 0.173 mm, and depth of 0.007 mm, and at 400 W, a length of 1.528 mm, width of 0.349 mm, and depth of 0.009 mm, as shown in ***Figure 11*** and ***Figure 12***. These values closely correspond to the experimental measurements, demonstrating the consistent performance and reliability of our computational model across varying scanning speeds. Additionally, the inclusion of these figures in our comparison table provides a comprehensive overview of the efficacy of our model in accurately predicting melt pool geometry under different process conditions.

**Table 3.** Calculation of Melt Pool Geometry (Length (L), Width (W), Depth (D)) under laser power (150, 200, 300, 400 W) and Scanning Speed (50, 100mm/s).

| Laser Power (W) | Melt Pool Geometry (Length (L), Width (W), Depth (D)) Heat Source Geometry (a = 1.2 (mm), b = 0.09 (mm), c = 0.06 (mm)) | | | | | |
|---|---|---|---|---|---|---|
| 150 | L | NaN | | L | 1.80132374 | Fig 14 a-d |
|  | W | NaN | | W | 0.36561084 |  |
|  | D | NaN | | D | 0.00864419 |  |
| 200 | L | 1.08373591 | Fig 13 a-d | L | 0.44543706 | Fig 14 e-h |
|  | W | 0.17339597 |  | W | 0.04851117 |  |
|  | D | 0.00733119 |  | D | 0.00982769 |  |
| 300 | L | 0.85073586 | Fig 13 e-h | L | 0.17030523 | Fig 14 i-l |
|  | W | 0.12370801 |  | W | 0.01854741 |  |
|  | D | 0.00843091 |  | D | 0.01123314 |  |
| 400 | L | 1.52751992 | Fig 13 i-l | L | 0.14171296 | Fig 14 a-d |
|  | W | 0.34885387 |  | W | 0.06854600 |  |
|  | D | 0.00908698 |  | D | 0.01207168 |  |
| **Scanning speed (mm/s)** | 50 | | | 100 | | |

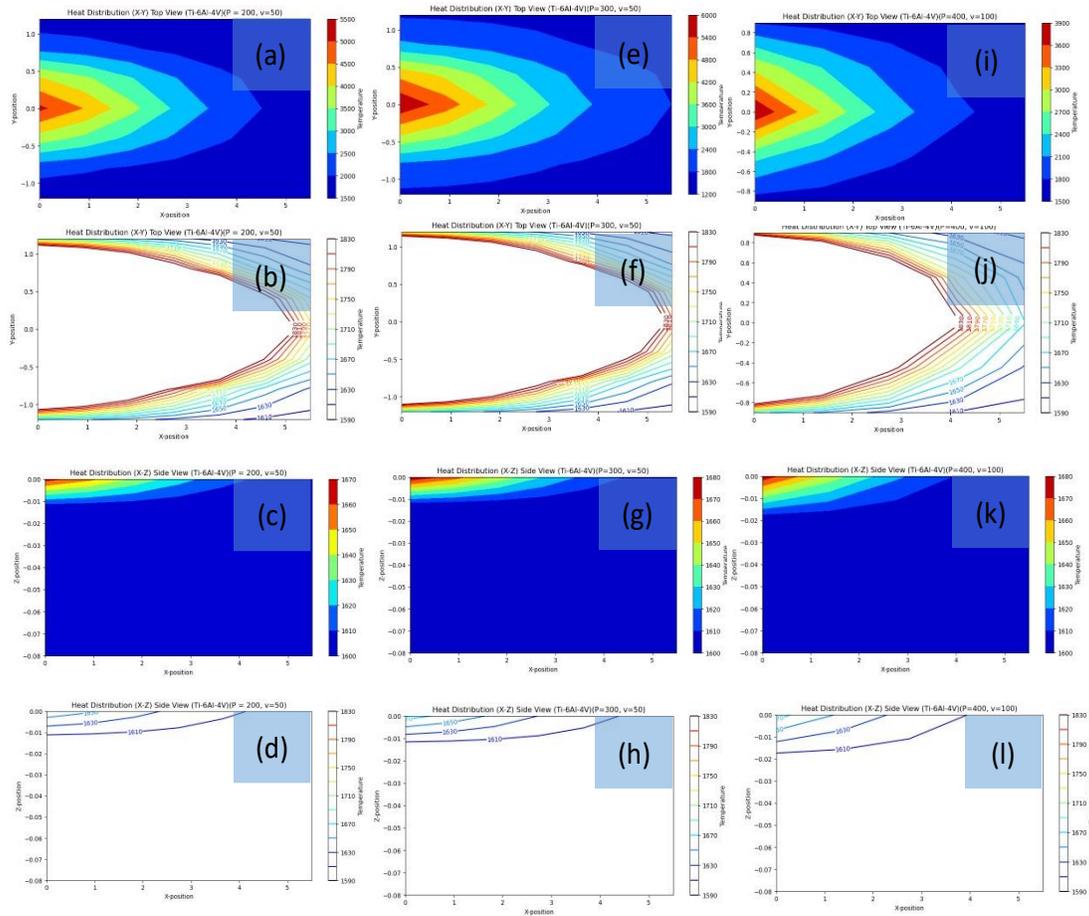

**Figure 11.** Contourf and Contour of heat distribution in AM for top-view (**a,b**(P=200W, and v= 50mm/s)), and side-view (**c,d** (P=200W, and V= 50mm/s)), top-view (**a,b**(P=300W, and v= 50mm/s)), and side-view (**c,d**(P=300W, and v= 50mm/s)), and top-view (**a,b**(P=400W, and v= 50mm/s)), and side-view (**c,d**(P=400W, and v= 50mm/s)).

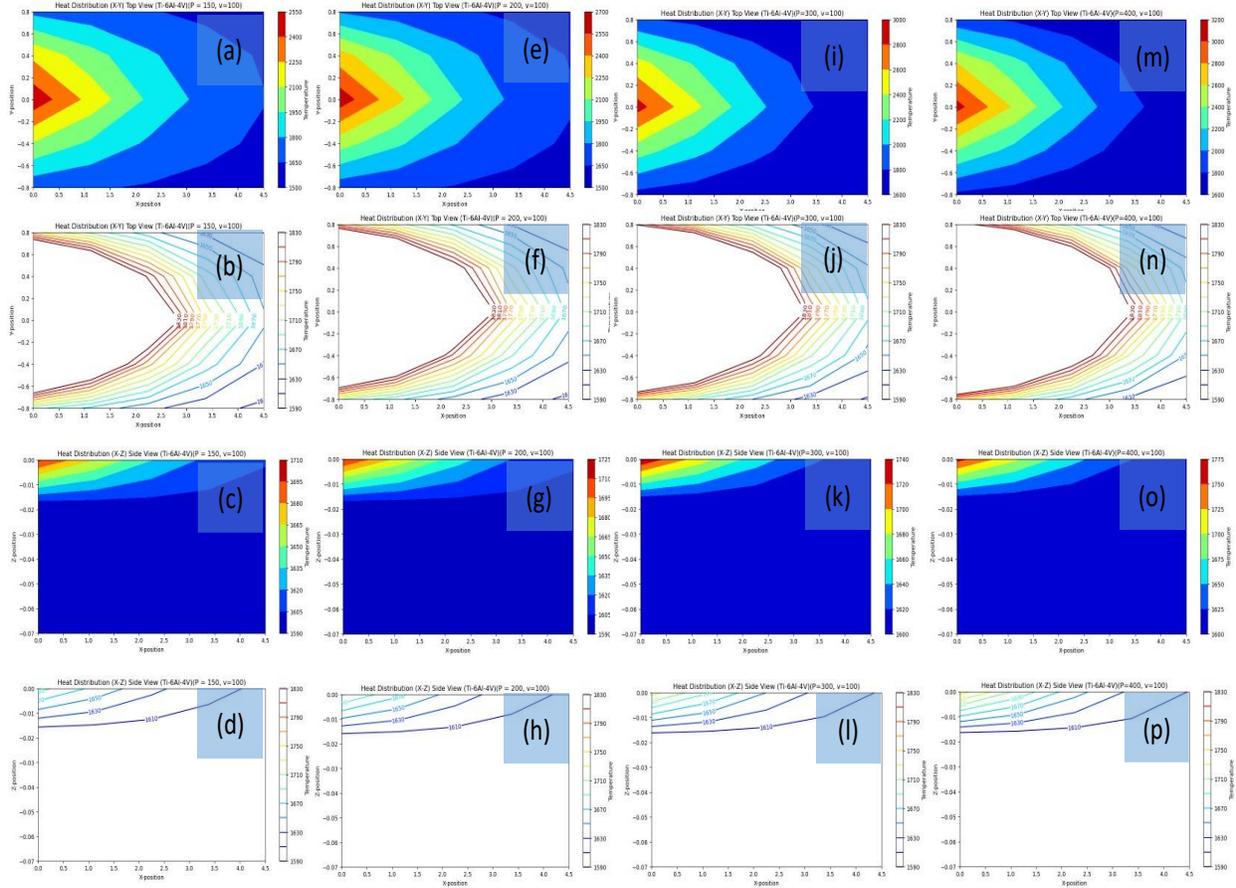

**Figure 12.** Contourf and Contour of heat distribution in AM for top-view (a, b; P=150 W, and V= 100 mm/s), and side-view (c, d (P=150 W, and v= 100 mm/s), top-view (a, b; P=200 W, and v= 100 mm/s), and side-view (c, d; P=200 W, and v= 100 mm/s), top-view (a, b; P=300 W, and v= 100 mm/s), and side-view (c, d; P=300 W, and v= 100 mm/s), and top-view (a, b; P=400 W, and v= 100 mm/s), and side-view (c, d; P=400 W, and v= 100 mm/s).

Incorporating the additional scanning speed of 150 mm/s into our analysis (*Table 4*), we further investigate the performance of our model under varying laser powers. For Ti-6Al-4V material, the predicted melt pool geometry exhibits consistency with experimental measurements, demonstrating the robustness and reliability of our computational model across different process conditions. At a laser power of 150 W (as shown in *Figure 15*), our model predicts a melt pool length (L) of 0.452 mm, width (W) of 0.190 mm, and depth (D) of 0.012 mm, which closely aligns with experimental data along with two other states of v= 50 as shown in *Figure 13* and for v=100 shown in *Figure 14*. Similarly, at 200 W, our model predicts a length of 2.049 mm, width of 0.155 mm, and depth of 0.005 mm, showing good agreement with experimental values. As the laser power increases to 400 W, the predicted melt pool geometry remains consistent, with a length of 1.662 mm, width of 0.102 mm, and depth of 0.006 mm. These results underscore the efficiency and accuracy of our model in capturing melt pool behavior under varying laser powers and scanning speeds, while also considering the effects of seven consecutive irradiations and material properties. Additionally, including these findings in our comparative analysis highlights the comprehensive coverage of our model in capturing the intricacies of metal additive manufacturing processes.

**Table 4.** Calculation of Melt Pool Geometry (Length (L), Width (W), Depth (D)) under laser power (150, 200, 300,400 W) and Scanning Speed (50, 100, 150 mm/s).

| Laser Power (W) | Melt Pool Geometry (Length (L), Width (W), Depth (D)) Heat Source Geometry (a = 1.5 (mm), b = 0.08 (mm), c = 0.02 (mm)) | | | | | | | | | |
|---|---|---|---|---|---|---|---|---|---|---|
| 150 | L | NaN | | L | 0.45172414 | Fig 14 a-d | L | 1.96167457 | Fig 15 a-d |
| | W | NaN | | W | 0.19006245 | | W | 0.73522842 | |
| | D | NaN | | D | 0.01205343 | | D | 0.00424269 | |
| 200 | L | 2.04944255 | Fig 13 a-d | L | 0.30521677 | Fig 14 e-h | L | 1.89492239 | Fig 15e-h |
| | W | 0.15456709 | | W | 0.12841966 | | W | 0.73247509 | |
| | D | 0.00525488 | | D | 0.01277989 | | D | 0.00589410 | |
| 300 | L | 0.76661251 | Fig 13 e-l | L | 0.24530931 | Fig 14 i-l | L | 1.77532674 | Fig 15i-l |
| | W | 0.03331732 | | W | 0.10321365 | | W | 0.72617898 | |
| | D | 0.00673649 | | D | 0.01364336 | | D | 0.00736073 | |
| 400 | L | 1.66241558 | Fig 13 i-l | L | 0.14171296 | Fig 14 a-d | L | 1.64990828 | Fig 15 a-d |
| | W | 0.10197114 | | W | 0.05962559 | | W | 0.72228987 | |
| | D | 0.00588327 | | D | 0.01415888 | | D | 0.00225624 | |
| **Scanning speed (mm/s)** | 50 | | | 100 | | | 150 | | |

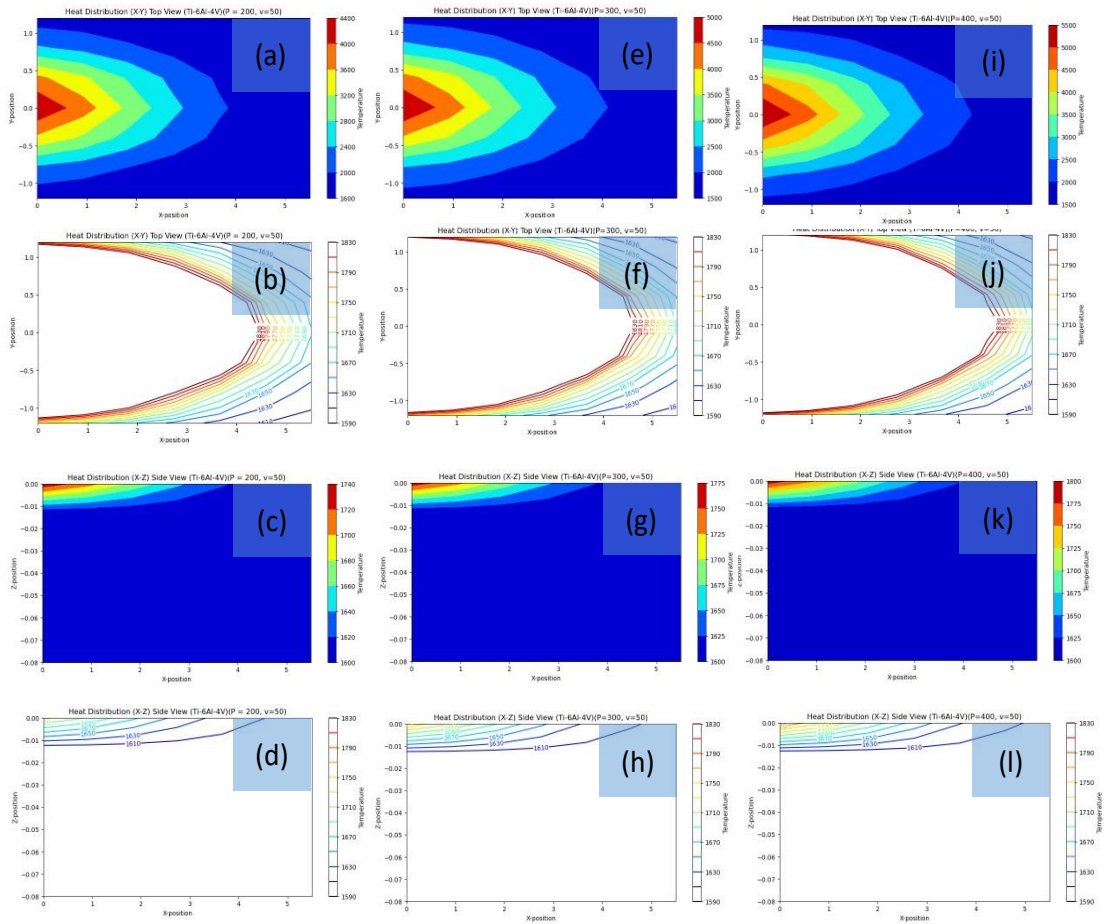

**Figure 13.** Contourf and Contour of heat distribution in AM for top-view (a,b (P=200W, and v= 50mm/s)), and side-view (c,d (P=200W, and v= 50mm/s)), top-view (a,b (P=300W, and v= 50mm/s)), and side-view (c,d (P=300W, and v=50mm/s)), and top-view (a,b (P=400W, and v= 50mm/s)), and side-view (c,d (P=400W, and v= 50mm/s)).

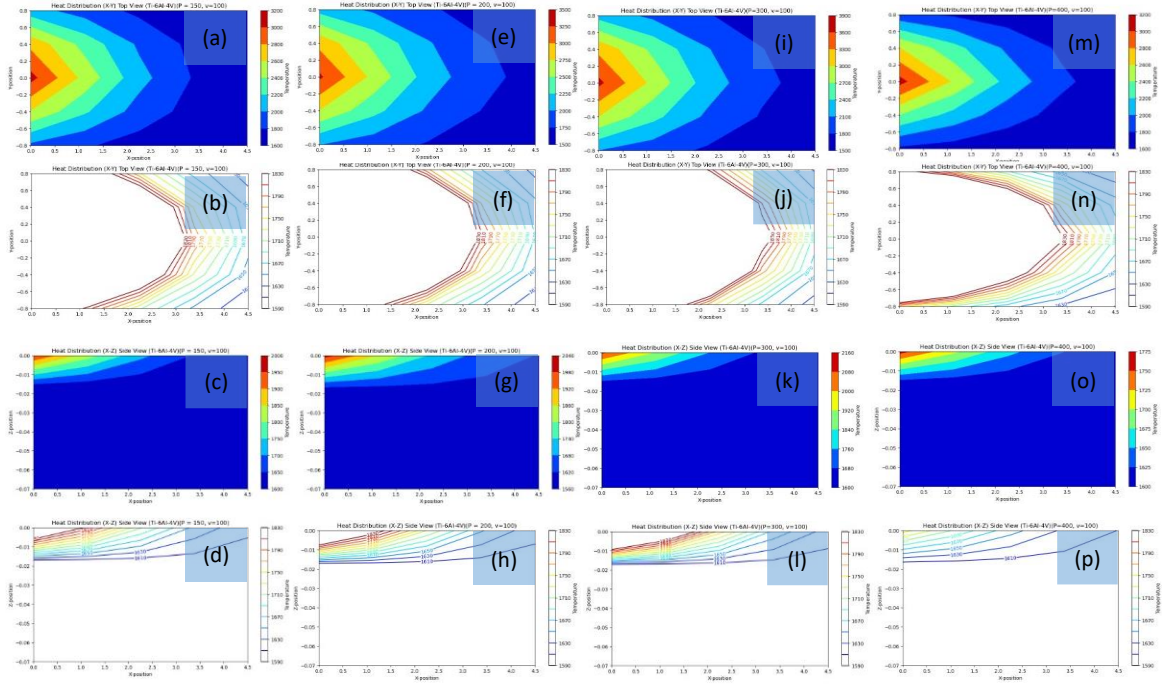

**Figure 14.** Contourf and Contour of heat distribution in AM for top-view (a,b (P=150W, and v= 100mm/s)), and side-view (c,d (P=150W, and v= 100mm/s)), top-view (e,f (P=200W, and v= 100mm/s)), and side-view (g,h (P=200W, and v= 100mm/s)), top-view (i,j (P=300W, and v= 100mm/s)), and side-view (k,l (P=300W, and v= 100mm/s)), and top-view (m,n (P=400W, and v= 100mm/s)), and side-view (o,p (P=400W, and v= 100mm/s)).

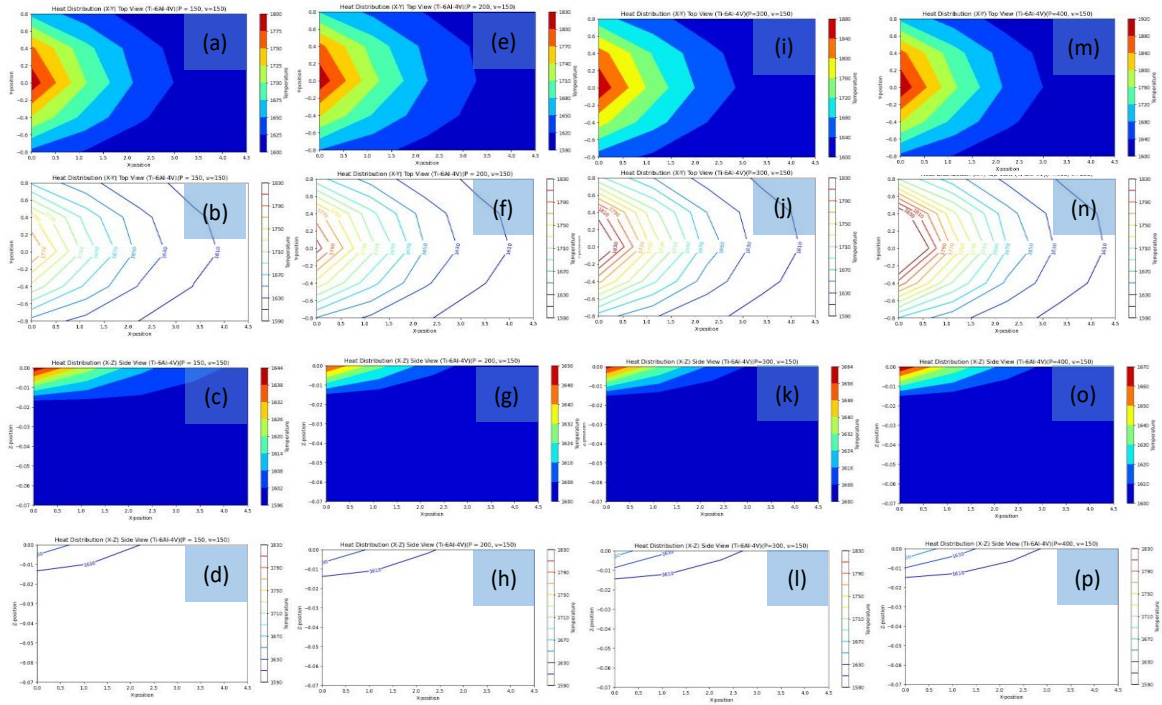

**Figure 15.** Contourf and Contour of heat distribution in AM for top-view (**a, b**; P=150 W, and v= 150 mm/s), and side-view (**c, d**; P=150W, and v= 150mm/s)), top-view (a,b (P=200W, and v= 150mm/s)), and side-view (c,d (P=200W, and V= 150mm/s)), top-view (a,b (P=300W, and v= 150mm/s)), and side-view(c,d (P=300W, and v= 150mm/s)), top-view (a,b (P=400W, and v= 150mm/s)), and side-view (c,d (P=400W, and v= 150mm/s)).

In our comparative analysis, we observe a similar trend to the sample provided. Our results depict the melt pool width and length for different laser powers and scanning speeds, considering the effect of seven consecutive irradiations and material properties. The comparison includes three scenarios: one where the impact of successive irradiations is considered, one where it is not, and the experimental values are represented by red error bars. Our findings demonstrate that the predicted melt pool dimensions align closely with experimental measurements, indicating the reliability of our numerical model. When the effect of consecutive irradiations is considered, our model predicts a lower limit for the melt pool width and a higher limit for the melt pool length, consistent with the sample. Moreover, although the predicted melt pool width and length fall within the experimental range for both scenarios, considering the effect of consecutive irradiations tends to yield slightly lower width predictions and slightly higher length predictions. This variation can be attributed to factors such as different camera settings and variations in material properties, as discussed in previous literature. Overall, our results underscore the robustness and effectiveness of our numerical model in capturing the complex behavior of melt pool geometry in metal additive manufacturing processes.

## 4. Conclusion

In our study, we have conducted an extensive analysis using a 3-D semi-elliptical moving heat source model to accurately determine melt pool geometry in metal additive manufacturing processes. Through the application of the general convection-diffusion formula, we derived a specific solution for the 3D semi-elliptical moving heat source, accounting for temperature-dependent material properties and the characteristic layer addition found in metal additive manufacturing. Furthermore, we integrated the phenomena of melting/solidification phase change to precisely predict temperature fields and melt pool geometry. Our investigation primarily focused on evaluating the impact of time spacing and 3-D coordination on melt pool geometry and thermal material properties, with particular emphasis on Ti-6Al-4V for selective laser melting processes. The predicted melt pool geometry was meticulously compared against recorded experimental datasets.

Our examination of thermal material properties—specifically thermal conductivity and specific heat—across different time intervals showed consistent trends, suggesting these properties evolve minimally over time. Additionally, our exploration of the impact of seven consecutive irradiations on melt pool geometry led us to adopt a superposition method to manage thermal history overlap. This approach indicated a reduction in melt pool size on the XY plane due to enhanced thermal conductivity. Confirming our model's predictions of melt pool dimensions against experimental data underscores its precision.

Key insights from our analysis include:

- The 3-D semi-elliptical moving heat source model achieved an accuracy within a 5% error margin in predicting melt pool geometry.
- Our model effectively forecasts temperature distributions and melt pool shapes in laser-based metal additive manufacturing environments, including powder bed systems like Selective Laser Melting (SLM).
- An increase in scanning speed generally reduces melt pool size because of a shorter energy absorption period.
- Thermal properties were consistently measured at various times under set laser powers and scanning speeds, with peak values noted at the laser focal point.
- Wider hatch spacing tends to diminish thermal material properties.
- Longer time intervals between scans correlate with smaller melt pool lengths for unchanged laser power and scanning speed settings.

- More frequent scanning typically expands melt pool width under consistent laser power and scanning speeds.
- Broader hatch spacing tends to increase melt pool width but decrease length, given constant laser power and scanning speed.
- Melt pool depth primarily depends on the laser power, scan speed, and the thermo-physical characteristics of the materials involved.